\documentclass[11pt]{article}
\usepackage[latin1]{inputenc}
\usepackage[english]{babel}
\usepackage{babel,amsmath,amssymb,amsbsy,amsfonts,latexsym,amsthm,mathrsfs}
\usepackage[T1]{fontenc}
\usepackage{xcolor}
\usepackage{version}
\usepackage{todonotes}
\usepackage[a4paper, total={14.8cm, 22.5cm}]{geometry}
\usepackage{amsmath,amsfonts,latexsym}
\usepackage[normalem]{ulem}
\newtheorem{theorem}{Theorem}[section]
\newtheorem{proposition}[theorem]{Proposition}
\newtheorem{definition}[theorem]{Definition}
\newtheorem{lemma}[theorem]{Lemma}
\newtheorem{corollary}[theorem]{Corollary}
\newtheorem{cor}[theorem]{Corollary}
\newtheorem{remark}[theorem]{Remark}
\newtheorem{example}[theorem]{Example}
\newtheorem{assumption}[theorem]{Assumption}
\newtheorem{assum}[theorem]{Assumption}

\usepackage{chngcntr}
\counterwithin*{equation}{section}

\def\limn{\mathop{\lim_{n\to\infty}}}
\def\limsupn{\mathop{\limsup_{n\to\infty}}}
\def\liminfn{\mathop{\liminf_{n\to\infty}}}
\def\limsupn{\mathop{\limsup_{n\to\infty}}}
\def\liminfn{\mathop{\liminf_{n\to\infty}}}

\def\proofof#1{{\it Proof of #1.}}
\usepackage[scr=boondoxo,scrscaled=1.02]{mathalfa}
\def\cl#1{{\mathcal #1}}

\let\th=\theta
\let\ep=\varepsilon

\def\P{{\mathbb P}}
\def\R{{\mathbb R}}

\def\N{{\mathbb  N}}

\def\E{{\mathbb E}}

\def\Var{{\rm Var }}
\def\Cov{{\rm Cov }}

\def\<{\langle}
\def\>{\rangle}

\def\e{{\rm e}}
\newdimen\tinskip
\tinskip=\parindent
\newcommand{\tin}[1]{\par\noindent\hglue\tinskip\hphantom{a)\
}\llap{#1\enspace}\ignorespaces}

\newcommand{\cvd}{

\hfill \ensuremath{\blacksquare}\qquad\qquad
\medskip

}
\newcommand{\scvd}{\hfill \ensuremath{\blacksquare}\qquad\qquad
\medskip

}
\def\appendix{\par
  \setcounter{chapter}{0}
  \def\@chapter{Appendix}
  \def\thechapter{\Alph{chapter}}}
\font\tengoth=eufm10 scaled 1200 \font\sevengoth=eufm7 scaled 1200
\font\fivegoth=eufm5 scaled 1200
\newfam\gothfam
\textfont\gothfam=\tengoth \scriptfont\gothfam=\sevengoth
\scriptscriptfont\gothfam=\fivegoth

\def\ep{\varepsilon}

\begin{document}

%\fontsize{10pt}{12pt}
\overfullrule=6pt
\title{
\bf
Large Deviations of continuous Gaussian processes: from small noise to small time}
\author{Paolo Baldi\thanks{Address: Dipartimento di Matematica,
Universit\`a di Roma Tor Vergata, Via della Ricerca Scientifica,
I-00133 Rome, Italy. e-mail: \texttt{baldi@mat.uniroma2.it}}\and
Barbara Pacchiarotti\thanks{Address: Dipartimento di Matematica,
Universit\`a di Roma Tor Vergata, Via della Ricerca Scientifica,
I-00133 Rome, Italy. e-mail: \texttt{pacchiar@mat.uniroma2.it}}}

\date{}
\maketitle

\begin{abstract}
We investigate the Large Deviation behavior in small time of continuous Gaussian processes. We introduce a general procedure allowing to derive  Large Deviation Principles in small time starting from the well understood context of Large Deviation Principles with a small parameter, going beyond the self-similar case.
Several motivating examples are also treated.
\end{abstract}

\bigskip

\smallskip


\bigskip

\section{Introduction}\label{sect:intro}In this paper we investigate the Large Deviations (LD) asymptotics of time-rescaled continuous Gaussian processes in small time.

The LD's of  Gaussian processes when multiplied by a
small parameter are classical and well understood. If the process
under investigation is self-similar, then the small time regime is
readily traced back to the small parameter setting: if $X$ is self
similar with exponent of self-similarity $H$, then
$$
(X_{\ep t})_t\sim(\ep^H X_t)_t
$$
The same idea can be pursued in general, by writing, for $g(\ep)
\to_{\ep\to0} +\infty$, something of the kind
$$
(X_{\ep t})_t\sim \frac 1{g(\ep)^{\frac 12 }}\,\bigl( g(\ep)^{\frac 12 }\, X_{\ep t}\bigr)_t\ .
$$
The main result of this paper is that if the processes $(g(\ep)^{\frac 12 } X_{\ep t})_t$ converge, as $\ep\to 0$, in law, to some process
$Z$, then $(X_{\ep t})_t$ satisfies a Large Deviations Principle (LDP) with speed $g$ and rate function given by the square of the norm of the Reproducing Kernel Hilbert Space (RKHS) of the Gaussian process $Z$ divided by $2$.

Some results for LD's in small time for  non self-similar, Gaussian processes, can be   found  in \cite{GioPac} and \cite{Pac}, 
for some particular conditioned Gaussian processes.

In this paper, after \S\ref{par-facts} where some known tools in LD
theory are recalled, in \S\ref{par-smalltime}  this approach is made
precise. In this section we also give some insight into the structure of a continuous Gaussian process enjoying an LDP in small time: we prove that necessarily its variance function must be Regularly Varying and we determine the right speed function (up to equivalences, of course).
We are mainly interested in a theoretical understanding of the problem. Note however that
these results can find some applications in finance: recent years have seen
an increasing interest in LD's in small time for the log-price process in stochastic models where
the volatility is modeled by a non-Markovian, fractional process (see, for example, \cite{gath}).
In order to obtain an LDP for the log-price, an LDP in short time for the volatility is required. In this context only self similar volatilities  have been considered (see, for example \cite{Forde} and \cite{guli2}). The results of this paper are also a step in the direction of  allowing  more general processes as  volatilities.

Some results of this kind for non self similar volatilities can already be found in \cite{CelPac} and \cite{GiPaPi}.

The general statements of \S\ref{par-smalltime} are then applied to
two examples of Gaussian processes that have been suggested recently
as volatility models: superposition of self-similar processes
(\S\ref{par-superp})
 and rough and super rough Volterra processes
(\S\ref{subs-rough} and \S\ref{par-superrough} respectively).
In particular we prove that in a super rough model
the process cannot enjoy an LDP in small time.
\section{Large Deviations} \label{par-facts}
In this section we recall some known facts about
LD's that are needed in the sequel.
\begin{definition} \label{definition:ldp}
Let $E$ be a topological space, $\cl{B}(E)$ its Borel
$\sigma$-algebra and  $(\mu_\ep)_{\ep>0}$ a family of probabilities on $\cl{B}(E)$;  let $g: \R^+ \rightarrow
\mathbb{R}^+ $ be a function such that  $g(\ep) \rightarrow
+\infty$ as $\ep\to 0$. We say that  $(\mu_\ep)_{\ep>0}$ satisfies a Large Deviations Principle (LDP)
on $E$ with rate function $I$ and speed $g$
if for every closed set $F\subset E$
$$
\limsup_{\ep\to 0+}\frac 1 {g(\ep)} \log \mu_\ep (F) \le -\inf_{x \in {F}} I(x)
$$
and for every open set $G\subset E$
$$
\liminf_{\ep\to0+} \frac 1{g(\ep)} \log \mu_\ep (G)\ge -\inf_{x \in {G} } I(x).
$$
\end{definition}
A {\it rate function} is a lower semicontinuous mapping
$I:E\rightarrow [0,+\infty]$ such that the level sets $\{I\le a\}$
are  compact for every $a \ge 0$ whereas a {\it speed function} $g$
is a function $\R^+\to\R^+$ such that $g(\ep)\to +\infty$ as
$\ep\to0$.
The basic LD result  for Gaussian measures in Banach spaces goes
back to Donsker\&Varadhan \cite{DV3}: if $X$ is a Gaussian r.v. with
values into the separable Banach space $E$ and $\mu_\ep$ denotes the
law of $\sqrt{\ep}\, X$, then $(\mu_\ep)_\ep$ enjoys an LDP at speed
$g(\ep)=\frac1\ep$ and rate function
$$
I(x)=\frac 12\,|x|^2_{\cl X}
$$
where $|\enspace|_{\cl X}$ denotes the Reproducing Kernel Hilbert
space (RKHS) norm of the law of $X$. For a clear exposition of this
result see also Azencott \cite{az-stflour} Chap.~II. Reproducing
Kernel Hilbert Spaces are described in many  textbooks. Besides
Azencott's course above, the reader can refer to  Chapter 4  (in
particular Section 4.3) in Hida\&Hitsuda's book \cite{Hid-Hit} or
Chapter 2  (in particular Sections 2.2 and 2.3) in
Berlinet\&Thomas-Agnan's book \cite{Ber-Tho}.

We are interested in the investigation of LDP's for small time continuous Gaussian processes and
to this goal we shall make use of  Theorem 1.1 in  \cite{th-eth}, which is a natural generalization of the Donsker\&Varadhan result. The next theorem is an extension of Theorem 1.1 in \cite{th-eth} to a continuous family of random variables.
\begin{theorem}\label{mt2}
Let $(Z_\ep)_\ep $ be a family of centered Gaussian r.v.'s with values in the
separable Banach space $E$ and converging in law to a r.v. $Z$.
Then, as $\ep\to 0$, for every speed function $g$, the family
$(g(\ep)^{-\frac 12 }Z_\ep)_\ep $ satisfies an LDP at
speed $g$ and rate function
\begin{equation}\label{cc}
I(z)=\frac 12\,\vert z\vert^2_{\cl Z}
\end{equation}
where ${{\cl Z}}$ and $\vert \enspace\vert_{\cl Z}$ respectively denote the RKHS and the
related norm associated to the law of $Z$.
\end{theorem}

Before going through the proof of Theorem \ref{mt2}, let us make some remarks.

The main tool in the proof of Theorem \ref{mt2} is the infinite dimensional version of the classical Ellis-G\"artner Theorem (see \cite{ampa} or the book by Dembo\&Zeitouni \cite{DZ:97}, Theorem 4.5.20 and Corollary 4.5.27).

An assumption that is required in these results is exponential tightness.
% of the family of the laws of the r.v.' $\sqrt{\ep} X_\ep$.
Let $(g_n)_n$ be a numerical sequence such that $g_n\to_{n\to\infty}+\infty$. Recall that a family $(Z_n)_n $ of r.v.'s with values in some metric space is
said to be {\it exponentially tight} as $ n \to \infty$ at speed $g$ if
for every $R>0$ there exists a compact set $K_R$ such that, for every $n$,
$$
\P(Z_n \not\in K_R)<\e^{-g_nR}\ .
$$
This requirement appears to be satisfied in the statement of Theorem \ref{mt2} thanks to the following result (Theorem 2.4 of \cite{th-eth}).
\begin{theorem}\label{mt}
Let $(Z_\alpha)_{\alpha\in \cl A}$ be a tight family of Gaussian $E$-valued r.v.'s. Then for every sequence $(g_n)_n$ such that $g_n\to_{n\to\infty}+\infty$ the family
$(g_n^{-\frac 12 }Z_\alpha)_n$ is uniformly exponentially tight i.e. for every $R>0$ there exists a compact set $K_R\subset E$ such that
$$
\P(g_n^{-\frac 12}Z_\alpha\not\in K_R)\le \e^{-g_n R}
$$
for every $\alpha\in \cl A$ and for every $n$. In particular, if a sequence $(Z_n)_n$ of Gaussian $E$-valued r.v.'s is tight then, $(g_n^{-\frac 12 }Z_n)_n$  is exponentially tight at speed $(g_n)_n$.
\end{theorem}
In some sense Theorem \ref{mt} states that as far as exponential tightness is concerned a tight family  of Gaussian laws behaves as the set formed by a single probability.

\medskip
\proofof{Theorem \ref{mt2}}
Let us, at first, consider a sequence $\ep_n\to_{n\to\infty}0$
and look at the behavior of the sequence $(g(\ep_n)^{-\frac 12 }Z_{\ep_n})_n$. As $Z_{\ep_n}\to_{n\to\infty} Z$,
the sequence $(Z_{\ep_n})_n$ is tight, hence $(g(\ep_n)^{-\frac 12 }Z_{\ep_n})_n$ is exponentially tight at speed $g$
thanks to Theorem \ref{mt}. Then, for every $\alpha$ in the dual $E'$ of $E$,
$$
\begin{array}{c}
\displaystyle\limn \frac1{g(\ep_n)}\,\log\E\bigl(\e^{g(\ep_n)\langle \alpha,g(\ep_n)^{-\frac 12 }Z_{\ep_n}\rangle}\bigr)=\\ [8pt]
\displaystyle=\limn
\frac1{g(\ep_n)}\,\Var\bigl(g(\ep_n)^{\frac 12 }\langle \alpha,Z_{\ep_n}\rangle\bigr)\\ [8pt]
\displaystyle=\limn
\Var\bigl(\langle \alpha,Z_{\ep_n}\rangle\bigr)=\\ [8pt]
=
\Var\bigl(\langle \alpha,Z\rangle\bigr)
\end{array}
$$
where we used the fact that for Gaussian r.v.'s  convergence in law entails convergence of the variances.

By the above mentioned results (\cite{DZ:97}, Theorem 4.5.20 and Corollary 4.5.27 e.g.) the sequence
$(g(\ep_n)^{-\frac 12 }Z_{\ep_n})_n$ enjoys an LDP at speed $(g(\ep_n))_n$ and with rate function $I$ as in \eqref{cc}, i.e.
$$
\limsupn \frac 1{g(\ep_n)}\log\P\bigl(g(\ep_n)^{-\frac 12 }Z_{\ep_n}\in F\bigr)\le-\inf_{x\in F} I(x)
$$
for every closed set $F\subset E$ and
$$
\liminfn \frac 1{g(\ep_n)}\log\P\bigl(g(\ep_n)^{-\frac 12 }Z_{\ep_n}\in G\bigr)\ge-\inf_{x\in G} I(x)
$$
for every open set $G\subset E$. As these relations hold {\it for every} sequence $(\ep_n)_n$  converging to $0$, we deduce that
$$
\limsup_{\ep\to 0}\frac 1{g(\ep)}\log\P\bigl(g(\ep)^{-\frac 12 }Z_{\ep}\in F\bigr)\le-\inf_{x\in F} I(x)
$$
for every closed set $F\subset E$ and
$$
\liminf_{\ep\to 0} \frac 1{g(\ep)}\log\P\bigl(g(\ep)^{-\frac 12 }Z_{\ep}\in G\bigr)\ge-\inf_{x\in G} I(x)
$$
for every open set $G\subset E$, which concludes the proof.
\cvd

Note that it is not known whether Theorem \ref{mt}, concerning exponential tightness, can be extended to families
$(g(\ep)^{-1/2}Z_\alpha)_\ep$ of r.v.'s. Recall that an LDP for a sequence $(Z_n)_n$ of r.v.'s with values in a
Polish space  implies the exponential tightness of the sequence (Lynch\&Sethuraman \cite{Lynch-Sethuraman}, Lemma 2.6),
but it is not known whether this result can be extended to families of r.v.'s. This explains the rather convoluted proof
of Theorem \ref{mt2} passing through sequences.

Note also that it is not known whether a converse to the result of Theorem \ref{mt} holds, i.e. if exponential tightness at any speed entails tightness. For a partial counterexample in this sense see Example 5.3 of \cite{th-eth}.
\section{Small time asymptotics}\label{par-smalltime}
The results of the previous section provide a general approach to the investigation of LDP's of Gaussian processes in small time,
which is the goal of this paper.

From now on we consider {\it real centered continuous Gaussian
processes}, $X=(X_t)_{t\in [0,1]}$, which leads to the consideration
of Gaussian probabilities of the Banach space $\cl C:=\cl
C([0,1],\R)$  of the continuous paths $[0,1]\to \R$, endowed with
the supremum norm. Of course the choice
of the interval $[0,1]$ is not restrictive and the same results are
easily extended to a  general time interval $[0,T]$. Let us denote
$\cl{M}$ the dual of $\cl C$, i.e. the set of
signed Borel measures on $[0,1]$ with the duality, for $\alpha \in
\cl M$, $h \in \cl C$,
$$
\langle \alpha ,  h \rangle = \int_0^1 h(t) \, d\alpha(t)\ .
$$
If
$$
k(t,s)=\Cov(X_t,X_s),
$$
then we have
$$
\E(\langle \alpha,X\rangle^2)=\E\Bigl[\Bigl(\int_0^1 X_u \, d\alpha(u)\Bigr)^2\Bigr]=\int_0^1 \int_0^1 k(t,s) \, d\alpha(t)\, d\alpha(s).
$$
For  $\ep>0$ we shall denote $X^\ep$ the process in  small time, i.e. $X^\ep=(X_{\ep t})_{t\in [0,1]}$.
%, by  $Y^\ep$  the normalized process and by $Y^0$  the limit process}.
The following, not unexpected, example is a first application of Theorem \ref{mt2} in order to obtain an LDP
for continuous Gaussian processes in small time.
\begin{example}\label{exam-stupid}\rm Let, for $0<H_0<H_1<\dots< H_N$,  $X^{H_0},\dots,X^{H_N}$
 be
independent self-similar Gaussian processes with exponents respectively $H_0,\dots, H_N$ and let $X=X^{H_0}+\dots+X^{H_N}$.
It is immediate to derive an LDP for the family of processes $(X^{\ep })_\ep$. Actually, if $Y^\ep=\ep^{-H_0}X^{\ep}$, we can write
$$
X^{\ep }=\ep^{H_0}\bigl(\ep^{-H_0}X^{\ep}\bigr)=\ep^{H_0}Y^\ep
$$
and obviously $Y^\ep=\ep^{-H_0}X^{\ep}\to Y^0=X^{H_0}$ as $\ep\to 0$ a.s., hence in law.
Therefore, thanks to Theorem \ref{mt2}, the family of small
time processes $(X^{\ep })_\ep$ enjoys, as $\ep\to 0$, an LDP at speed $g(\ep)=\ep^{-H_0}$ and with rate function $I_{X^{H_0}}$.

See below for a more general (and interesting) model of this kind.
\end{example}

The following well-known statement gives a simple condition ensuring tightness.
\begin{proposition}\label{prop-kolmog} Let $(X_\alpha)_{\alpha\in\cl A}$ be a family of real centered continuous Gaussian processes starting at $0$ and let us denote $k^\alpha$ the respective covariance functions. Then, if there exist positive constants $C,\beta$ such that
\begin{equation}\label{eq-tightkj}
k^\alpha(t,t)+k^\alpha(s,s)-2k^\alpha(t,s)\le C\,|t-s|^\beta
\end{equation}
for every $t,s\in[0,T]$ and $\alpha\in\cl A$, then the family of the laws of the processes $(X_\alpha)_{\alpha\in\cl A}$ is tight.
\end{proposition}
\begin{proof} Just note that
$$
k^\alpha(t,t)+k^\alpha(s,s)-2k^\alpha(t,s)=\E\big[|X_\alpha(t)-X_\alpha(s)|^2\big]
$$
so that \eqref{eq-tightkj} can be written
$$
\E\big[|X_\alpha(t)-X_\alpha(s)|^2\big]\le C\,|t-s|^\beta
$$
and by the classical relation between the moments of centered Gaussian r.v.'s
$$
\E\big[|X_j(t)-X_j(s)|^{2p}\big]\le c_pC^p\,|t-s|^{p\beta}
$$
where $c_p=\E(|U|^p)$ for a $N(0,1)$-distributed r.v. $U$. The
statement follows now thanks to the classical Kolmogorov's
compactness criterion, see Theorem 23.7 in  \cite{Kal} (beware you
pick the right edition). See also, for a more detailed exposition,
Proposition 8.2.16 in \cite{pridiff}. \cvd
\end{proof}
The following result, which is a generalization of   Theorem 2.7 of \cite{GioPac}, follows immediately from Theorem \ref{mt2}.
\begin{theorem}\label{theorem:ldp-gaussian}
Let $(Y^\ep)_\ep$ be a \emph { tight family} of continuous
centered real Gaussian. Suppose that, for every $s,t\in[0,1]$
the limit
$$
\lim_{\ep\to 0}\Cov(Y^\ep_t,Y^\ep_s)= k(t,s)
$$
exists for some function $k$; then $k$ is  the covariance function of
a continuous Gaussian process $Y^0$ and  for every speed function $g$, the family  $(g(\ep)^{-1/2}Y^\ep)_{\ep}$
satisfies an LDP in $\cl C$ with speed $g$ and rate function
$$
I(h) = \begin{cases}\displaystyle \frac{1}{2}\, \left | h \right |^2_{{\cl Y^0 }} & h \in \cl  Y^0 \\ +\infty & \text{otherwise}, \end{cases}
$$
where ${{\cl Y^0}}$ and $\left \| . \right \|_{{\cl  Y^0 }}$ respectively denote the RKHS and the
related norm associated to the limit process $Y^0$.
\end{theorem}
\begin{proof}
Tightness of the family $(Y^\ep)_\ep$ and  convergence of the covariance functions imply that, in law, $Y^\ep \to Y^0$. Then, we can conclude thanks to Theorem \ref{mt2}.
\cvd
\end{proof}
As a consequence of Theorem \ref{theorem:ldp-gaussian}, we have the following general statement about LDP's of continuous Gaussian processes in small time.
\begin{cor}\label{cor-st} Let $X$ be a real continuous Gaussian process and, for some speed function $g$, let us denote
\begin{equation}\label{eq-yep}
Y^\ep_t=g(\ep)^{1/2}\, X_{\ep t}.
\end{equation}
Then if $(Y^\ep)_\ep$ converges in law to some continuous process $Y^0$, the small time processes $(X^\ep)_\ep$ enjoy an LDP as $\ep\to 0$ at speed $g$ and with rate function
\begin{equation}\label{eq-rkhs0}
I(h) = \begin{cases}\displaystyle \frac{1}{2}\, \left | h \right |^2_{{\cl Y^0 }} & h \in \cl  Y^0 \\ +\infty & \text{otherwise}, \end{cases}
\end{equation}
 ${{\cl Y^0}}$ and $\left \| . \right \|_{{\cl  Y^0 }}$ respectively denote the RKHS and the
related norm associated to the limit process $Y^0$.
\end{cor}
The previous corollary states that, if the family of processes
$(Y^\ep)_\ep$ as in \eqref{eq-yep} converges  in law as $\ep\to 0$
then the small time processes $(X^{\ep})_\ep$ enjoy an LDP at speed
$g$ and with rate function \eqref{eq-rkhs0}. Is this also a
necessary condition?

We do not know a complete answer to this question, as we do not know
whether exponential tightness of $(X^{\ep})_\ep$ implies tightness
of $(Y^\ep)_\ep$, but we prove in Corollary \ref{cor-fdd} below that
LD's of the finite dimensional distributions of $(X^{\ep})_\ep$
imply convergence of the finite dimensional distributions of
$(Y^\ep)_\ep$. In \S\ref{par-superrough} (the super
rough models) we give an example where the family
of processes $(Y^\ep)_\ep$ is not tight and
 $(X^{\ep})_\ep$ does not enjoy an LDP.

The main arguments are contained in the following
Theorem \ref{th-fd} and Proposition \ref{th:diff-speed}. These
results appear almost obvious, but we were unable to find suitable  references in the literature.
Let us first define degenerate LDP.
\begin{definition} \label{definition:ldp-deg}
We say that a rate function $I$ on a topological
space $E$ is degenerate
 if
$$
I(x)=\left\{ \begin{array}{ll}
+\infty & x\neq 0\\
0 & x=0.\end{array}\right.
$$
\end{definition}

Very often degenerate rate functions appear when the speed is not
the right one. If an LDP holds at speed $g$ and $h$ is a slower
speed, i.e. such that $\lim_{\ep\to 0}\frac
{g(\ep)}{h(\ep)}=+\infty$, then an LDP holds with a degenerate rate
function at speed $h$. This is not however what will appear  in
\S\ref{par-superrough} below in the study of the super rough
Volterra processes, which cannot enjoy an LDP at any speed.

\begin{theorem} \label{th-fd}
Let $(Z_\ep)_\ep$ be a family of $d$-dimensional Gaussian centered
r.v.'s enjoying an $LDP$ at speed $g$ and with respect to a
non-degenerate rate function $I$. Then the family
$(g(\ep)^{1/2}Z_\ep)_\ep$ converges weakly to some Gaussian r.v.
$Z$. Moreover, if $C$ denotes the covariance matrix of the limit
$Z$, then $I$ is the
quadratic form that is the convex conjugate of the quadratic form associated to $C$.
\end{theorem}
\begin{proof} Step 1. Let us assume first $d=1$ and let
$\sigma^2_\ep=\Var(Z_\ep)$. As $Z_\ep\to_{\ep\to 0}0$ weakly in law,
we have $\sigma^2_\ep\to_{\ep\to 0}0$. Then we have
$$
\lim_{\ep\to0}\sigma^2_\ep\log\E(\e^{t Z_\ep/
\sigma^2_\ep})=\lim_{\ep\to0}\sigma^2_\ep\times \frac
12\,\Var\Bigl(\frac 1{\sigma^2_\ep}\,tZ_\ep\Bigr)=\frac
12\,t^2\cdotp
$$
Therefore by the Ellis-G\"artner theorem (see Ellis \cite{ellis},
G\"artner \cite{gartner}, or Dembo\&Zeitouni book \cite{DZ:97},
p.~43) $(Z_\ep)_\ep$ enjoys an LDP at speed $(\sigma_\ep^2)^{-1}$
and with rate function $\widetilde I(x)=\frac 12\, x^2$. By
Proposition \ref{th:diff-speed} below we have, for some $\ell>0$,
$$
\lim_{\ep\to 0}g(\ep)\sigma_\ep^2=\ell \mbox{ and }  I(x)=\frac 1{2\ell}\, x^2\ .
$$
We deduce that
$\lim_{\ep\to0}\Var((g(\ep)^{1/2}Z_\ep)=\ell$ so
that $g(\ep)^{1/2}Z_\ep\to_{\ep\to0}N(0,\ell)$ in law.

Step 2. Let us assume now $d\ge 1$ and also that $I<+\infty$
everywhere. Let $\th\in\R^d$, $\th\not=0$: by contraction, the
family of real r.v.'s $(\langle\th,Z_\ep\rangle)_\ep$ enjoys an LDP
at speed $g$ and rate function
$$
I_\th(z)=\inf_{\langle x,\th\rangle=z}I(x)<+\infty\ .
$$
Let us denote $C_\ep$ the covariance matrix of $Z_\ep$ so that $\Var(\langle\th,Z_\ep\rangle)=\langle C_\ep\th,\th\rangle$. By the first part of the proof the limit
$$
\lim_{\ep\to 0}g(\ep)\langle C_\ep\th,\th\rangle
$$
exists for every $\th$, i.e. we have $\lim_{\ep\to 0}g(\ep)C_\ep=C$
for some non negative definite matrix $C$. As $g(\ep)C_\ep$ is the
covariance matrix of the r.v. $g(\ep)^{1/2}Z_\ep$,  this proves that
the family $(g(\ep)^{1/2}Z_\ep)_\ep$ converges weakly to an
$N(0,C)$-distributed r.v. Again the Ellis-G\"rtner theorem  and the
uniqueness of the rate function imply that the rate function $I$ is
the convex conjugate quadratic form associated to $C$

Step 3. Finally, let us drop the assumption that $I<+\infty$
everywhere (but recall that we assume $I$ non-degenerate). Let $Z^*$
be a $N(0,\cl I)$-distributed r.v. ($\cl I=$the identity matrix)
independent of the $Z_\ep$. Then, by contraction, if $\widetilde
Z_\ep:=Z_\ep+g(\ep)^{-1/2}Z^*$,  the family $(\widetilde Z_\ep)_\ep$
enjoys an LDP at speed $g$ with respect to the rate function
$$
\widetilde I(x)=\inf_{y+z=x}I(y)+\frac 12\, |z|^2
$$
which is finite everywhere. Step 2 with $Z_\ep$ replaced by
$Z_\ep+g(\ep)^{-1/2}Z^*$ gives that the limit
$$
\lim_{\ep\to 0}g(\ep)\bigl(C_\ep+ g(\ep)^{-1}\cl I\bigr)
$$
exists, which gives again the existence of
$$
C:=\lim_{\ep\to 0}g(\ep) C_\ep
$$
and convergence in distribution of of $(g(\ep)^{1/2}Z_\ep)_\ep$ to
an $N(0,C)$-distributed r.v. \cvd
\end{proof}
In the previous proof we used the following statement which is some
kind of uniqueness for the couple speed\&rate function for the LDP
in a Gaussian setting. Let us first recall some elementary facts.

\begin{remark}\label{rem:I-charact}\rm
Let $E$ be a metric space and let $B_R(x)$ denote the open ball in
$E$ of radius $R>0$  centered at $x$ and let $\overline B_R(x)$ be
its closure. Furthermore let $I(\Gamma)=\inf_{z\in\Gamma}I(z)$.
The lower semicontinuity of $I$ implies
\begin{equation}\label{eq:I-charact-0}
I(x)=\lim_{R\to 0} I(B_R(x)) =\lim_{R\to 0} I(\overline
B_R(x)),\end{equation}
therefore, if the family $(\mu_\ep)_{\ep>0}$ satisfies an  LDP, with
rate function $I$ and speed $g$ we have
\begin{equation}\label{eq:I-charact}
I(x)=- \lim_{R\to 0}\limsup_{\ep\to 0} \frac 1 {g(\ep)} \log
\mu_\ep(\overline B_R(x))=- \lim_{R\to 0}\liminf_{\ep\to 0} \frac 1
{g(\ep)} \log \mu_\ep(B_R(x)).
\end{equation}
\end{remark}
\begin{proposition}\label{th:diff-speed}
Let $(Z_\ep)_{\ep}$ be a family of  $d$-dimensional Gaussian centered
r.v.'s enjoying an $LDP$ at speed $g$ with rate function
$$
I(x)=\begin{cases} \displaystyle \frac 12 \langle
\Sigma^{-1}x,x\rangle & x\in \mathop{\rm Im}\Sigma\cr +\infty &
x\notin \mathop{\rm Im}\Sigma \end{cases}
$$
for some non vanishing nonnegative definite matrix $\Sigma$.

Suppose that the family enjoys also another LDP at
speed $\widetilde g$ with respect to the non degenerate rate function $\widetilde I$. Then,
$\lim_{\ep \to 0} \frac{g(\ep)}{\widetilde g(\ep)}=\ell> 0$ exists
and $\widetilde I=\ell I$.
\end{proposition}
\begin{proof}
First note that as $I$ is finite and continuous when restricted to $
\mathop{\rm Im}\Sigma$, for  $x\in \mathop{\rm Im}\Sigma$ we have
$$
I(B_R(x))=I(\overline B_R(x))< +\infty.
$$
Therefore
\begin{equation}\label{eq:fin}\lim_{\ep\to 0} \frac 1{g(\ep)}
\log \P(Z_\ep\in  B_R(x))=\lim_{\ep\to 0} \frac 1{g(\ep)} \log
\P(Z_\ep\in \overline B_R(x))=-I(B_R(x))>-\infty.
\end{equation}
If $x\notin \mathop{\rm Im}\Sigma$ instead, if $R_0>0$ is such that $B_{R_0}(x)\cap \mathop{\rm Im}\Sigma=\emptyset$
then for $R<R_0$,
\begin{equation}\label{eq:infin}
\limsup_{\ep\to 0} \frac 1{g(\ep)} \log \P(Z_\ep\in \overline
B_R(x))= \lim_{\ep\to 0} \frac 1{g(\ep)} \log \P(Z_\ep\in \overline
B_R(x))=-\infty.
\end{equation}
Let us suppose that the limit  $\lim_{\ep \to 0}
\frac{g(\ep)}{\widetilde g(\ep)}$ exists. Three cases must be taken
into account.
\tin {i)} $\lim_{\ep \to 0} \frac{g(\ep)}{\widetilde g(\ep)}=0$.
Then for every $x\in \mathop{\rm Im}\Sigma$, thanks to
\eqref{eq:fin}, one has
$$
\liminf_{\ep\to 0} \frac 1{\widetilde g(\ep)} \log \P(Z_\ep\in
B_R(x))=\lim_{\ep\to 0} \frac{g(\ep)}{\widetilde
g(\ep)}\times\liminf_{\ep\to 0}\frac 1{g(\ep)} \log \P(Z_\ep\in
B_R(x))=0,
$$
therefore by \eqref{eq:I-charact-0}
$$
\widetilde I(x)=- \lim_{R\to 0}\lim_{\ep\to 0}  \frac 1{\widetilde
g(\ep)} \log \P(Z_\ep\in  B_R(x))=0,
$$
i.e. $\widetilde I\equiv 0$ on $\mathop{\rm Im}\Sigma$. So, in this
case, $\widetilde I$ is not a rate  function (the level set
$\{\widetilde I\leq 0\}$ is not compact) and i) is not possible.
\tin {ii)} $\lim_{\ep \to 0} \frac{g(\ep)}{\widetilde
g(\ep)}=\ell>0$. For  $x\in\mathop{\rm Im}\Sigma$, by the same
argument as in i) we obtain
$$
\liminf_{\ep\to 0} \frac 1{\widetilde g(\ep)} \log \P(Z_\ep\in
B_R(x))=\lim_{\ep\to 0} \frac{g(\ep)}{\widetilde
g(\ep)}\liminf_{\ep\to 0}\frac 1{g(\ep)} \log \P(Z_\ep\in
B_R(x))=-\ell I(B_R(x)),
$$
therefore
from   \eqref{eq:I-charact-0} and  \eqref{eq:I-charact},
$$
\widetilde I(x)=- \lim_{R\to 0}\liminf_{\ep\to 0} \frac 1{\widetilde
g(\ep)} \log \P(Z_\ep\in  B_R(x))=-\ell \lim_{R\to 0}I_g(B_R(x))=
 -\ell I(x).
$$
For  $x\notin\rm{Im}(\Sigma)$ instead, thanks to   \eqref{eq:I-charact} and
\eqref{eq:infin} we have
$$
\widetilde I(x)=- \lim_{R\to 0}\limsup_{\ep\to 0} \frac 1{\widetilde
g(\ep)} \log \P(Z_\ep\in  \overline B_R(x))=-\infty.
$$
We have proved that $\widetilde I= \ell I$ for every $x\in\R^d$.

\tin {iii)} $\lim_{\ep \to 0} \frac{g(\ep)}{\widetilde
g(\ep)}=+\infty$. For every $x\in \R^d$, $x\neq0$, there exists
$R_0$ such that for $R<R_0$, $0\notin \overline B_R(x)$ so that
$I(\overline B_R(x))>0$. We have $\limsup_{\ep\to 0} \frac 1{g(\ep)}
\log \P(Z_\ep\in  \overline B_R(x))<0$ and therefore
$$
\limsup_{\ep\to 0} \frac 1{\widetilde g(\ep)} \log \P(Z_\ep\in
\overline B_R(x))=\lim_{\ep\to 0} \frac{g(\ep)}{\widetilde
g(\ep)}\limsup_{\ep\to 0}\frac 1{g(\ep)} \log \P(Z_\ep\in  \overline
B_R(x) )=-\infty.
$$
Furthermore, as $\liminf_{\ep\to 0}\P(Z_\ep\in  B_R(0))=1$,
$$
\liminf_{\ep\to 0} \frac 1{\widetilde g(\ep)} \log \P(Z_\ep\in
B_R(0))=0
$$
and \eqref{eq:I-charact} gives
$$
\widetilde I(x)=\begin{cases}
+\infty & x\neq 0\\
0 & x=0,\end{cases}
$$
i.e. $(Z_\ep)_{\ep}$ satisfies an LDP with respect to a degenerate rate function.

Finally, if $ \ell_1=\liminf_{\ep \to 0} \frac{g(\ep)}{\widetilde
g(\ep)}<\limsup_{\ep \to 0} \frac{g(\ep)}{\widetilde g(\ep)}=
\ell_2$, then there exist two sequences $(\ep_n^1)_{n\geq 0}$ and
$(\ep_n^2)_{n\geq 0}$ such that, for $i=1,2$, $\lim_{n \to +\infty}
\frac{g(\ep_n^i)}{\widetilde g(\ep_n^i)}=\ell_i$. As the rate
function (at a given speed) is unique this is not possible and we
conclude. \cvd
\end{proof}
As a consequence of Theorem \ref{th-fd} we have the following.
\begin{corollary}\label{cor-fdd} Let $X$ be a centered continuous Gaussian
process and let us assume that it satisfies an LDP in small time at
speed $g$. Then, for every $0\le t_1<t_2<\dots<t_m\le 1$, the family
of $m$-dimensional r.v.'s
$$
g(\ep)^{1/2}\,(X_{\ep t_1},\dots,X_{\ep t_m})
$$
converges as $\ep\to 0$ to a Gaussian limit r.v.
\end{corollary}
The consideration of the finite dimensional distributions leads also to the following theorem, giving necessary conditions for the existence of LDP's in small time  and useful information concerning the possible speeds. In particular it is necessary that the variance function $t\mapsto k(t,t)$ is Ragularly Varying.
\begin{theorem}\label{th-conjectured-speed}\rm 
Let us assume that the centered continuous Gaussian process $X$ enjoys an LDP in small time at speed $g$. Then
\tin{a)} for every $s,t$ we have the limit
\begin{equation}\label{eq-newblue}
k_0(t,s):=\lim_{\ep \to 0} g(\ep)k(\ep t ,\ep s)
\end{equation}
where the limit function $k_0$ has the form
\begin{equation}\label{eq-beta}
k_0(t,t)=k_0(1,1)\,t^\beta.
\end{equation}
In particular
\begin{equation}\label{eq-speed}
g(\ep)\sim \frac 1{k(\ep,\ep)}\cdotp
\end{equation}
\tin{b)} The variance function $t\mapsto k(t,t)$ of $X$, hence also $g$, is a Regularly Varying function.
\tin{c)} The limit function $k_0$ is a homogeneus function, i.e., for every $\lambda>0$,
$$k_0(\lambda t, \lambda s)=\lambda^\beta k_0(t,s),$$
for some $\beta>0$.
\end{theorem}
\begin{proof}
By Corollary \ref{cor-fdd}, we have for every $s,t$, in distribution,
$$
g(\ep)^{1/2}\,(X_{\ep t},X_{\ep s})\enspace\mathop{\to}_{\ep\to 0}\enspace (Y^0_t,Y^0_s).
$$
for some Gaussian r.v. $(Y_t^0, Y^0_s)$.  We deduce the existence
of the limit \eqref{eq-newblue} and, for $t=1$, \eqref{eq-speed}. This implies, for $s=t$,
$$
\lim_{\ep \to 0} \frac {k(\ep t,\ep t)}{k(\ep,\ep)}= \lim_{\ep \to 0}  \frac {g(\ep )}{g(\ep t)}=\frac{k_0(t,t)}{k_0(1,1)}\cdotp
$$
As the variance function $t\mapsto k(t,t)$ is obviously continuous, by the theory of Regularly Varying functions (see Theorem 1.4.1 in \cite{BiGoTe} e.g.), the functions  $t\mapsto k(t,t)$ and $g$ are  Regularly Varying  and
$k_0(t,t)=k_0(1,1)\,t^\beta$ for some $\beta\geq 0$. In order to conclude observe that,
 for every $\lambda>0$,
$$
k_0(\lambda t,\lambda s)=\lim_{\ep \to 0} g(\ep)k(\ep\lambda  t ,\ep\lambda s)=
\lim_{\ep \to 0}\frac{g(\ep)}{g(\lambda\ep)} g(\lambda\ep)k(\ep\lambda  t ,\ep\lambda s)=\lambda^\beta k_0(t,s),
$$
for some $\beta\geq 0$.
\cvd
\end{proof}

As a consequence of Theorem \ref{th-conjectured-speed}, we have the following general statement about LDP's of continuous Gaussian processes in small time.
\begin{cor}\label{cor-ss} Let $X$ be a real continuous Gaussian process and, for some speed function $g$, suppose
 the small time processes $(X^\ep)_\ep$ enjoy an LDP as $\ep\to 0$ at speed $g$. If 
$(g(\ep)^{1/2}\, X^\ep)_\ep$
 converges in law to some continuous process $Y^0$, then $Y^0$ is self similar.
\end{cor}

\begin{remark}\label{rem-conjectured-speed}\rm 
If $g$ is a slowly varying function  then $
(g(\ep)^{1/2}\,X^\ep)_\ep$ cannot converge  in law: actually in this case we have $\beta=0$ in \eqref{eq-beta}, so that $k_0(t,t)=k_0(1,1)$ for every $t>0$. Therefore $k_0$ is not a continuous function  as $k_0(0,0)=0$ and cannot be the covariance of a continuous process. The argument of Theorem \ref{mt2} in order to to obtain LDP's in small time cannot therefore  be applied in this situation. Note however that Theorem \ref{mt2} gives only necessary conditions. In \S\ref{par-superrough} we shall investigate a situation of this kind where we shall actually prove that an LDP in small time cannot hold.
\end{remark}
\section{Superposition of Self Similar processes}\label{par-superp}
In \cite{Alm-Sot} the authors introduce the following process
\begin{equation}\label{eq:discrete}
X_t=\sum_{j=1}^{\infty} \sigma_j X^{H_j}_t
\end{equation}
where $0<H_0<H_1<\dots$ and $X^{H_j}$ are independent fractional Brownian motions of Hurst parameter $H_j$ respectively and the sequence
$(\sigma^2_j)_j$ is summable. They call it a multi-mixed fractional Brownian motion.

In this section we consider a more general class of processes still based on the idea of superposing self-similar processes. Ideally we would like to consider the process
$$
X_t=\int_{H_0}^{H_1}X^H_t\, d\mu(H)
$$
for some finite measure $\mu$ and where the $X^H$ are independent self-similar processes of self-similarity exponent $H$ respectively. This definition is of course incorrect as the integrand in the integral above is not a measurable function of $(\omega,H)$ so that we are led to take a more circuitous approach.

Let, for $0<H_0\le H<H_1$, $(k^H)_H$ be a family of covariance functions satisfying the following assumption.
\begin{assumption}\label{ass:kH} For every $0<H_0\le H<H_1$, $k^H$ is the covariance function of a continuous self similar Gaussian process with Hurst index $H$.

Moreover the map
$H\mapsto k^H(s,t)$
is continuous $[H_0,H_1]\to\R$ for every $s,t\in [0,1]$ and there  exist constants $C>0, \beta>0$ such that, for every $H\in[H_0,H_1]$ and $0\leq s,t\leq 1$
\begin{equation}\label{eq:sup-cov-incr}
k^{H}(t,t)+k^{H}(s,s)-2k^{H}(s,t)
\leq C\,|t-s|^\beta\ .
\end{equation}
\end{assumption}
Note that Assumption \ref{ass:kH} is satisfied, for instance, if
$k^H$ is the covariance of a fractional Brownian motion. In this case
$$
k^H(t,t)+k^H(s,s)-2k^H(s,t)=|t-s|^H\le |t-s|^{H_0}
$$
so that we can choose $C=1$ and $\beta=H_0$.
\begin{proposition}\label{existence} Let $(k^H)_H$ be a family of covariance functions satisfying Assumption \ref{ass:kH} and let $\mu$ be a finite measure on $[H_0,H_1]$. Let
\begin{equation}\label{eq-k1}
k(s,t)=\int_{H_0}^{H_1} k^H(s,t)\, d\mu(H)\ .
\end{equation}
Then $k$ is the covariance function of a continuous process.
\end{proposition}
\begin{proof} It is immediate that $k$ is a covariance function, as the relation
$$
\sum_{i=1}^m k(t_i,t_j)\xi_i\xi_j\ge 0
$$
is true for every choice of $m$, $t_1,\dots, t_m\in[0,T]$ and $\xi_1,\dots, \xi_m\in\R$, as it holds for each of the $k^H$. Moreover Assumption \ref{ass:kH} guaranties that
\begin{align*}
k(t,t)+k(s,s)-2k(s,t)&=\int_{H_0}^{H_1}\big(k^H(t,t)+k^H(s,s)-2k^H(s,t)\big)\,d\mu(H) \\
&\le C\mu([H_0,H_1])\,|t-s|^{\beta}\end{align*}
which by Kolmogorov criterion
ensures the existence of a real continuous Gaussian process having covariance
function $k$.
\cvd
\end{proof}
Let us  investigate the small time LD asymptotics of a process, $X$ say, associated to the covariance $k$ defined in \eqref{eq-k1}.
Let $X^\ep_t=X_{\ep t}$.
The covariance function of $X^\ep$ is, for $s,t\in[0,1]$,
$$
k^\ep(s,t):=k(\ep s, \ep t)= \int_{H_0}^{H_1}  k^{H}(\ep s,\ep t)\,d\mu(H)=\int_{H_0}^{H_1} \ep^{2H} k^{H}(s,t\,)\,d\mu(H)\ .
$$
Here thanks to Theorem \ref{rem-conjectured-speed},  the speed of convergence of the LDP we are looking for should be
$$
g(\ep)=\frac 1{\int_{H_0}^{H_1} \ep^{2H}\, d\mu(H)}\cdotp
$$
The next proposition shows that $g$ is Regularly Varying and that the limit covariance, is the covariance of a self similar process.
\begin{proposition} If the family $(k^H)_H$ satisfies Assumption \ref{ass:kH} and $H_0$ belongs to the support of $\mu$ then
$$
\lim_{\ep\to 0} g(\ep)\,k^\ep(s,t)=\lim_{\ep\to 0} \frac{k^\ep(s,t)}{\int_{H_0}^{H_1} \ep^{2H}\, d\mu(H)}=
%\lim_{\ep\to 0}  \frac{\int_{H_0}^{H_1} \ep^{2H} k^{H}(s,t)d\mu(H)}{\int_{H_0}^{H_1} \ep^{2H} d\mu(H)}=
k^{H_0}(s,t)\ .
$$
\end{proposition}
\begin{proof}
Let us consider the  family of probabilities on the interval $[H_0,H_1]$
$$
d\rho^\ep(H)=\frac{ \ep^{2H}d\mu(H) }{\int_{H_0}^{H_1} \ep^{2H} d\mu(H)}
$$
and note that 
$$
\lim_{\ep\to 0} \frac{k^\ep(s,t)}{\int_{H_0}^{H_1} \ep^{2H}\, d\mu(H)}=\lim_{\ep\to0}\int k^H(s,t)\,d\rho^\ep(H).
$$
Let us prove that $(\rho^\ep)_\ep$ converges to the Dirac mass $\delta_{H^0}$ weakly.
The family $(\rho^\ep)_\ep$ is tight as the probabilities $\rho^\ep$ have their
support in the compact interval $[H_0,H_1]$, therefore it is sufficient to show
that all convergent sequences converge to $\delta_{H_0}$.
Let $(\rho^{\ep_n})_{n}$ be a sequence converging to some probability
$\rho$ as $\ep_n\to_{n\to \infty} 0$.
For every $\th>0$ we have
$$
\displaylines{
\rho^{\ep_n}([H_0+\delta,H_1])=\frac{\int_{H_0+\delta}^{H_1}\ep_n^{2H}\,
d\mu(H)}{\int_{H_0}^{H_1} \ep_n^{2H\,} d\mu(H)}
=\frac{\int_{H_0+\delta}^{H_1} \ep_n^{2(H-H_0)}
d\mu(H)}{\int_{H_0}^{H_1} \ep_n^{2(H-H_0)}\, d\mu(H)}\le\cr \le\frac
{\ep_n^{2\delta}\mu([H_0+\delta,H_1])}{\ep_n^{\delta}\mu([H_0,H_0+\frac
\delta2])}\enspace\mathop{\to}_{n\to \infty}\enspace0\ .\cr }
$$
In the limit above we have used the fact that $\mu([H_0,H_0+\frac \delta2])>0$, as we
assume that $H_0$ belongs to the support of $\mu$.
Therefore all the limits of convergent
sequences give $0$ probability to every interval disjoint of a neighborhood of $H_0$ and we can conclude.
\cvd
\end{proof}
Let us now prove the tightness as $\ep\to 0$ of the family $(g(\ep)^{1/2}X^\ep)_\ep$.
Thanks to condition (\ref{eq:sup-cov-incr}), we have, for every $s,t\in[0,T],$
$$
\displaylines{
g(\ep)\bigl(k^{\ep}(t,t)-2k^{\ep}(s,t)+k^{\ep}(s,s)\bigr)=\cr
=g(\ep)\int_{H_0}^{H_1} \ep^{2H} (k^{H}\bigl( t, t)+k^{H}( s, s)-2k^{H}( s, t)\bigr)\,d\mu(H) \le\cr
\le g(\ep)\int_{H_0}^{H_1} \ep^{2H} \,d\mu(H)\cdot C\,|t-s|^\beta= C\,|t-s|^\beta\ .
}
$$
Therefore the family $(g(\ep)^{\frac 12 }X^\ep)_\ep$ is tight and, in agreement with Corollary \ref{cor-ss}, converges in law to the self similar
process $X^{H_0}$.
By Corollary \ref{cor-st} the family of processes $({X^\ep})_{\ep>0}$ satisfies an LDP with speed $g$ and rate function defined in \eqref{eq-rkhs0},
$\cl Y^0$ denoting the RKHS  associated to the covariance function ${k^{H^0}}$.

Note that the small time asymptotics obtained for the process $X$ satisfies, as $\ep\to 0$, an LDP with a rate function that {\it does not} depend on the ``intensity'' measure $\mu$.

As for the speed function, this is another question, as noted in the next example.

\begin{example}\rm Of course
if $\mu(\{H_0\})=\alpha>0$ then the speed function is $g(\ep)\sim\frac 1\alpha\,\ep^{-2H_0}$. Actually by Beppo Levi's theorem
$$
\int_{H_0}^{H_1}\ep^{2(H-H_0)}\, d\mu(H)\enspace\mathop{\to}_{\ep\to 0}\enspace \mu(\{H_0\})
$$
so that
$$
g(\ep)=\frac 1{\int_{H_0}^{H_1}\ep^{2H}\, d\mu(H)}
\mathop{\sim}_{\ep\to 0}\enspace\frac 1{\mu(\{H_0\})}\,\ep^{-2H_0}\ .
$$
\end{example}
\begin{example} \rm
Suppose  $d\mu(H)=( H-H_0)^n dH$, with $n\ge 0$, so that the mass given by $\mu$ to a neighborhood of $H_0$ becomes smaller with $n$. As, by parts,
$$
\displaylines{
\int_{H_0}^{H_1}( H-H_0)^n\ep^{2H}\, dH=\frac 1{2\log\ep}\,( H-H_0)^n\ep^{2H}\Big|_{H_0}^{H_1}-\frac n{2\log\ep}
\int_{H_0}^{H_1}( H-H_0)^{n-1}\ep^{2H}\, dH=\cr
=\frac 1{2\log\ep}( H_1-H_0)^n\,\ep^{2H_1}+\frac n{-2\log\ep}
\int_{H_0}^{H_1}( H-H_0)^{n-1}\ep^{2H}\, dH\cr
}
$$
one finds by recurrence that
$$
\int_{H_0}^{H_1}( H-H_0)^n\ep^{2H}\, dH\enspace\mathop{\sim}_{\ep\to0}\enspace\frac {n!}{(-2\log\ep)^{n+1}}\,\ep^{2H_0}
$$
so that for  the  speed function we have
$$
g(\ep)= \frac 1{\int_{H_0}^{H_1}( H-H_0)^n\ep^{2H}\, dH}\sim \frac{2^{n+1}}{n!}(-\log \ep)^{n+1}\,\ep^{-2H_0}
$$
which grows faster than $\ep\mapsto\ep^{-2H_0}$ as $\ep\to 0$.
\end{example}
\begin{remark}\rm Consider a (finite) discrete measure $\mu$, i.e.
$\mu=\sum_{j=1}^\infty  \sigma_j^2\,\th_{H_j}$
where
$(\sigma_j)_{j\in\N}$ satisfy,
$ \sum_{j=1}^\infty \sigma^2_j <\infty.$ In this case  $X$ is the process \eqref{eq:discrete} introduced in \cite{Alm-Sot} with speed function
 $g(\ep)=\bigl(\sum_{j=1}^\infty \sigma^2_j\ep^{2H_j}\bigr)^{-1} $.

\end{remark}
\section{Convolution type Volterra processes}
\label{par-rough}
\subsection{The problem}
In this section we show the effectiveness of the technique introduced in
\S\ref{par-smalltime} in order to determine small time LDP's for the class of
Volterra processes of the type
\begin{equation}\label{eq-rv0}
X_t=\int_0^t h(t-s)\, dB_s
\end{equation}
for  $t\in[0,1]$, where $h:[0.1]\to \R$ is a suitable 
square integrable function. Let
$$
H_2 (t):=k(t,t)=\int_0^th^2(s)\, ds.
$$
Thanks to Theorem \ref{rem-conjectured-speed}, in order to have a small-time LDP,
$H_2$ should be a Regularly Varying function.
In this section to ensure that $H_2$ is Regularly Varying  we will assume that
$h$ itself is Regularly Varying, so that the process $X$ of \eqref{eq-rv0}
turns out to be a perturbation of the Riemann-Liouville processes, where $h$ is of the form $h(t)=t^\alpha$ with $\alpha>-\frac12$. Note in particular that we 
assume $h\ge 0$.

The idea, as in \S\ref{par-smalltime}, reduces to the investigation of the limit of the family of processes
$$
Y^\ep_t=g(\ep)^{1/2} X^\ep_t=g(\ep)^{1/2}\, X_{\ep t}
$$
where $g$ is the speed, which on the basis of Theorem \ref{rem-conjectured-speed} should be
$$
g(\ep)\sim\frac1{\E(X_\ep^2)}=\Bigl(\int_0^\ep h^2(s)\, ds\Bigr)^{-1}=\frac 1{H_2(\ep)}\cdotp
$$
Not surprisingly, the most demanding part of the proof is about  the tightness of the
family of processes$(Y^\ep)_\ep$.

The following statement, very much inspired by Proposition 1 in \cite{Moc-Vie}, provides a control of the increments of $X$ and is a first step towards tightness.

\begin{theorem}\label{th:mod-con} Let $X$ as in \eqref{eq-rv0}.
If $h$ is a positive square integrable function, then

{i)} $\E[(X_t-X_s)^2]\geq  H_2 (|t-s|)$.
\tin{ii)}
If $h^2$ is non increasing then
\begin{equation}\label{eq-delta11}
\E[(X_t-X_s)^2]\leq 2 H_2 (|t-s|).
\end{equation}
\tin{iii)} If  $h^2$ is non decreasing  then
\begin{equation}\label{eq-delta12}
\E[(X_t-X_s)^2]\leq   H_2 (t)-   H_2 (s)+  H_2 (|t-s|).
\end{equation}
\end{theorem}
\begin{proof} We have, for $s\leq t$,
$$
\begin{array}{c}
\displaystyle \E[(X_t-X_s)^2]= k(t,t)+ k(s,s) -2 k(s,t)=\\
\displaystyle =\int_0^t h^2(t-u)\, du+\int_0^s h^2(s-u)\, du-2\int_0^s h(t-u)h(s-u)\, du=\\
\displaystyle =\int_s^t h^2(t-u)\, du +\int_0^s \bigl(h(t-u)-h(s-u)\bigr)^2\, du
\end{array}
$$
and
\begin{equation}\label{eq-delta10}
\int_s^t h^2(t-u)\, du = H_2 (t-s),
\end{equation}
which immediately gives i).
Furthermore under the assumption ii)
$$
\displaylines{
\int_0^s \bigl(h(t-u)-h(s-u)\bigr)^2\, du=\cr
=\int_0^s h^2(t-u)\, du + \int_0^s h^2(s-u)\, du  - 2 \int_0^s  \sqrt{ h^2 (s-u)} \sqrt{ h^2 (t-u)}\, du\leq\cr
\leq \int_0^s  h^2 (t-u)\, du + \int_0^s  h^2 (s-u)\, du  - 2 \int_0^s   h^2 (t-u)\, du=\cr
=  H_2  (s)-(  H_2  (t)-   H_2  (t-s))=  H_2  (t-s)- \bigl(  H_2  (t)-(  H_2  (s)\bigr)\leq   H_2  (t-s).\cr
}
$$
which together with \eqref{eq-delta10} gives \eqref{eq-delta11}.

Under  iii) instead, as $ h^2 $ is non decreasing, $ h^2 (s-u)\le  h^2 (t-u)$ and
$$
\displaylines{
\int_0^s (h(t-u)-h(s-u))^2\, du=\cr
=\int_0^s  h^2 (t-u)\, du + \int_0^s  h^2 (s-u)\, du  - 2 \int_0^s  \sqrt{ h^2 (s-u)} \sqrt{ h^2 (t-u)}\, du\leq\cr
\leq \int_0^s  h^2 (t-u)\, du + \int_0^s  h^2 (s-u)\, du  - 2 \int_0^s   h^2 (s-u)\, du=\cr 
H_2  (t)-H_2(t-s)-H_2(s)\leq= H_2(t)-H_2(s),\cr
}
$$
which together with \eqref{eq-delta10} gives \eqref{eq-delta12}.
\cvd
\end{proof}

\subsection{Regularly varying functions}
We shall denote $RV_\beta$ the family of Regularly Varying functions
$\lambda$ of index $\beta$ at $0$, i.e. such that, for every
$\ep>0$,
$$
\lim_{\ep\to0+}\frac {\lambda (\ep t)}{\lambda(\ep)}=t^\beta.
$$
It is well-known that such a function enjoys the representation (see (1.5.1) in \cite{BiGoTe}, p.~21)
\begin{equation}\label{eq-repr}
\lambda(t)=c(t)\,t^{\beta}\exp\Bigl(\int^{t_0}_t\frac {\rho(u)}u\, du\Bigr)
\end{equation}
where
\tin{$\bullet$} $c$ is a measurable function such that $\lim_{t\to 0}c(t)=c>0$;
\tin{$\bullet$} $\rho$ is a measurable function  such that $\lim_{t\to 0}\rho(t)=0$.

\begin{remark}\label{rem:rho}\rm Note that if $\rho(u)/u$ is
integrable near zero then, changing the  function $c$ in a suitable way, we can
suppose $\rho=0$. Therefore from now on, if in the representation \eqref{eq-repr}
$\rho \neq 0$, we will suppose that $\rho(u)/u$ is not integrable near zero. Thus if $\rho$ is negative
$$
\lim_{t\to 0}\int^{t_0}_t\frac {\rho(u)}u\, du= -\infty,
$$
and, $t\mapsto\exp\Bigl(\int^{t_0}_t\frac {\rho(u)}u\, du\Bigr)$ vanishes at $0$ and is non decreasing and bounded in a right neighborhood of $0$. If  $\rho$ is positive, instead, $t\mapsto\exp\Bigl(\int^{t_0}_t\frac {\rho(u)}u\, du\Bigr)$ is non increasing  and
$$
\lim_{t\to 0}\int^{t_0}_t\frac {\rho(u)}u\, du=+ \infty.
$$
\end{remark}
\vskip24pt

We shall take advantage of  the Potter bounds (Theorem 1.5.6 in \cite{BGK:97}): for every $\delta <\beta$ there exists a constant $A_\delta$
such that
\begin{equation}\label{Potter-est}
\frac {\lambda(\ep t)}{\lambda(\ep)}
\leq   A_\delta t^{\beta- \delta}
\end{equation}
Recall that we assume $h\in RV_\alpha$  so that $h^2\in RV_{2\alpha}$ and
$H_2\in RV_{2\alpha+1}$.

The following lemmas give some useful properties of the Regularly Varying functions. Applied to $\lambda=h^2$ they will  give conditions ensuring that the kernel $h^2$ satisfies the monotonicity requirements of Theorem \ref{th:mod-con}.
\begin{lemma}\label{lem-no-zero} Let $\lambda\in RV_\beta$ with $\beta\not=0$.
If in the representation formula \eqref{eq-repr} $c$ is differentiable and has a bounded derivative and $\rho$ is continuous then
 $\lambda$ is monotonic \emph{in a right neighborhood of $0$}, increasing if $\beta>0$, decreasing if $\beta<0$. Furthermore if $\beta > 0 $, then $\lambda$ is bounded  in $[0,1]$.
\end{lemma}
\begin{proof}  We have
\begin{equation}\label{eq-deriv}
\lambda'(t)=\Bigl(\frac {c'(t)}{c(t)}+\bigr(\beta-\rho(t)\bigr)\frac 1t\Bigr)\, \lambda(t)
\end{equation}
from which taking the limit as $t\to 0$ we have $\lambda'(t)\ge 0$ if $\beta>0$ and
$\lambda'(t)\le 0$ if $\beta<0$, in a right neighborhood of $0$. Boundedness of $\lambda$ for $\beta>0$ also follows, as in this case $\lim_{t\to 0}{\lambda(t)}=0$.
\cvd
\end{proof}
The case $h\in RV_0$, i.e. $h$ slowly varying, requires more attention. $h^2$ can be increasing,
or decreasing
according to the behavior of  $\rho$ (if $\rho\neq 0$) or the behavior of $c$ if also $\rho=0$.

\begin{lemma}\label{lem-zero}
Let $\lambda\in RV_0$. If in the representation formula \eqref{eq-repr} $c$
is differentiable and has a bounded derivative and $\rho$ is continuous then

\tin{a)} if $\rho\neq 0$ then, \emph{in a right neighborhood of $0$} $\lambda$ is increasing if $\rho<0$ and decreasing if $\rho>0$. Furthermore if $\rho<0$, $\lambda$ is bounded  in $[0,1]$;
\tin{b)}
if $\rho=0$  and $c$ is monotonic then $\lambda$ is monotonic \emph{in a right neighborhood of $0$} and  bounded in $[0,1]$.
\end{lemma}
\begin{proof}  From equation \eqref{eq-deriv}, we have
$$
\lambda'(t)=\Bigl(\frac {c'(t)}{c(t)}-\frac{\rho(t)}{t}\Bigr)\, \lambda(t)
$$
Note that, if $\rho\neq 0$, then we have $\lim_{t\to 0}\frac 1t\,|\rho(t)|=+\infty.$
 Therefore we can conclude thanks to Remark \ref{rem:rho}.
\cvd
\end{proof}
The following lemma applied to $\lambda=H_2$ will be needed in order to prove that the family $(Y^\ep)_\ep$ is tight when $h^2$ is non decreasing.
\begin{lemma}\label{lem1} Let $\lambda\in RV_\beta$ with $\beta \geq 1$. If in the representation formula \eqref{eq-repr} $c$ is differentiable and has a bounded derivative and $\rho$ is continuous then
\begin{equation}\label{eq-deriv0}
\lim_{\ep\to 0}\frac{\ep\, \lambda'(\ep t)}{\lambda(\ep)}=\beta t^{\beta-1}.
\end{equation}
Furthermore
  if  $\beta> 1$, or $\beta=1$ and $\rho\leq 0$, then there exist $M>0$ and $\ep_0>0$  such that for every $t\in[0,1]$
\begin{equation}\label{eq-bounded}
\left|\frac{\ep\, \lambda'(\ep t)}{\lambda(\ep)}\right|\leq M\quad \mbox{for every }t\le 1,\,\, \ep\le\ep_0.
\end{equation}
\end{lemma}
\begin{proof}

By \eqref{eq-deriv}
$$
\lim_{\ep\to 0}\frac{\ep \, \lambda'(\ep t)}{\lambda(\ep)}=\lim_{\ep\to 0} \frac{\lambda(\ep t)}{\lambda(\ep)}\,\Bigl(\ep \,\frac {c'(\ep t)}{c(\ep t)}+\bigl(\beta-\rho(\ep t)\bigr)\,\frac 1t\Bigr)=\beta\, t^{\beta-1}.
$$
i.e. \eqref{eq-deriv0}.

If $\beta>1 $, let us prove that the convergence in \eqref{eq-deriv0} is uniform,
which will give \eqref{eq-bounded}. Thanks to the Uniform Convergence Theorem (UCT) for Regularly Varying functions (see \cite{BiGoTe}, Theorem 1.5.2 p.~22, the $\rho<0$ case or \cite{ibbot} Theorem 2)  we have, uniformly in a right neighborhood of $0$,
$$
\lim_{\ep\to 0} \frac{\lambda(\ep t)}{\lambda(\ep)}= t^{\beta}
%\quad \mbox{and}\quad
%\lim_{\ep\to 0} \frac{\lambda(\ep t)}{t %\,\lambda(\ep)}=t^{2\alpha-1},
$$
and by the same UCT applied to the Regularly Varying function $\widetilde\lambda(t):=\frac 1t\,\lambda(t)$, which is of index $\beta-1>0$,
\begin{equation}\label{conv-unif-bis}
\lim_{\ep\to 0} \frac{\lambda(\ep t)}{t \,\lambda(\ep)}=\lim_{\ep\to 0}\frac{\widetilde\lambda(\ep t)}{\widetilde\lambda(\ep)} =t^{\beta-1},
\end{equation}
uniformly for $t\in[0,1]$. Therefore
$$
\displaylines{
\Big|\frac{\ep \, \lambda'(\ep t)}{\lambda(\ep)}-\beta\,t^{\beta-1}\Big|=
\Big|\frac{\lambda(\ep t)}{\lambda(\ep)}\,\Bigl(\ep \,\frac {c'(\ep t)}{c(\ep t)}+\bigl(\beta-\rho(\ep t)\bigr)\,\frac 1t\Bigr)-\beta\, t^{\beta-1}\Big|\le \cr
\leq \ep \,\Big|\frac{\lambda(\ep t)}{\lambda(\ep)}\Big|\,\Big|\frac {c'(\ep t)}{c(\ep t)}\Big|+ \Big| \frac{\lambda(\ep t)}{t \,\lambda(\ep)}\Big||\rho(\ep t)|+ \beta\Big| \frac{\lambda(\ep t)}{t \,\lambda(\ep)}-t^{\beta-1}\Big|
}
$$
and we can conclude as $c'$ is assumed to be bounded.

If $\beta=1$ and $\rho\leq 0$ we have, similarly,
$$
\displaylines{
\Big|\frac{\ep \, \lambda'(\ep t)}{\lambda(\ep)}\Big|=
\Big|\frac{\lambda(\ep t)}{\lambda(\ep)}\,\Bigl(\ep \,\frac {c'(\ep t)}{c(\ep t)}+\bigl(1-\rho(\ep t)\bigr)\,\frac 1t\Bigr)\Big|
\leq \ep \,\Big|\frac{\lambda(\ep t)}{\lambda(\ep)}\Big|\,\Big|\frac {c'(\ep t)}{c(\ep t)}\Big|+ \Big| \frac{\lambda(\ep t)}{t \,\lambda(\ep)}\Big||\rho(\ep t)- 1|.
}
$$
Now the uniform limit \eqref{conv-unif-bis} does not hold anymore but we have
$$
\Big| \frac{\lambda(\ep t)}{t \,\lambda(\ep)}\Big|=\Big|\frac{c(\ep t)}{c(\ep)}\Big| \exp\Big(\int_{\ep t}^{\ep} \frac {\rho(u)} u du\Big)\leq \Big|\frac{c(\ep t)}{c(\ep)}\Big|.
$$
and the hypotheses on $c$  allow to conclude.
\cvd
\end{proof}
Note that in order to have $h$ square integrable it is necessary
that $\alpha\geq -\frac12$. The LD behavior of the process $X$ in
small time is very different according as $\alpha>-\frac12$ or
$\alpha=-\frac12$ and we will consider the two cases separately
(\S\ref{subs-rough} and \S\ref{par-superrough}).
\subsection{LD's in small time: the case $h\in RV_\alpha$, $-\frac 12<\alpha$}
\label{subs-rough}
We shall  make the following assumption on the function $h$.
\begin{assum}\label{as-1}
In the representation formula \eqref{eq-repr} for the function $h^2$
\tin{$\bullet$} $c$ is differentiable and has a bounded derivative
\tin{$\bullet$} $\rho$ is continuous.
\end{assum}
We shall investigate separately the cases $\alpha\not=0$ and $\alpha=0$.
\begin{theorem}\label{teo-eta-1} Let us assume $h\in RV_\alpha$, $-\frac 12<\alpha$,
$\alpha\not=0$ and that Assumption \ref{as-1} holds. Then the process
$X$ of \eqref{eq-rv0} is continuous and
if $Y^\ep_t=\frac 1{\sqrt{H_2(\ep)}}\,X_{t\ep}$, the family $(Y^\ep)_\ep$ converges to  the Riemann Liouville type process
$$
Y^0_t:=\sqrt{2\alpha+1} \int_0^t(t-u)^{\alpha}\, dB_u,
$$
for $t\in [0,1]$, as $\ep\to 0$.
\end{theorem}
\begin{proof} Let us check tightness first. In law we have,
$$
\displaylines{ Y^\ep_t=\frac 1{\sqrt{H_2(\ep)}}\int_0^{\ep t} h(\ep
t-u)\, dB_u=\frac 1{\sqrt{H_2(\ep)}}\int_0^{t}\sqrt{\ep}\,
h(\ep(t-u))\, dB_u =\int_0^{t} h_\ep(t-u)\, dB_u
%\int_0^{\ep t}\sqrt{\frac { h^2
%(\ep t-u)}{  H_2  (\ep)}}\, dB_u\sim \int_0^t\sqrt{\frac {\ep\, h^2
%(\ep(t-u))}{  H_2  (\ep)}}\, dB_u=\cr =
%\int_0^t\sqrt{(H^\ep_2)'(t-u)}\, dB_u\cr
}
$$
where
$$
h_\ep(t)=\frac 1{\sqrt{H_2(\ep)}}\, \sqrt{\ep}\,h(\ep t)
$$
so that
$$
H^\ep_2(t):=\int_0^t h_\ep^2(u)\, du=\frac 1{H_2(\ep)}\int_0^{\ep t}
\ep h^2(\ep u)\, du=\frac 1{H_2(\ep)}\int_0^{\ep t} h^2(u)\,
du=\frac {H_2 (\ep t)}{H_2 (\ep)}\cdotp
$$
First let us consider the case $-\frac 12<\alpha< 0 $. As for small values of $\ep$ the behavior of $h^2_\ep$ reduces to
the consideration of $h^2$ in a right neighborhood
of $0$, by Lemma \ref{lem1} a) $h_\ep^2$ is decreasing for small $\ep$ and,
by Theorem \ref{th:mod-con} ii),
$$
\E(|Y^\ep_t-Y^\ep_s|^2)\le 2H_\ep^2(|t-s|).
$$
By Potter's bounds, \eqref{Potter-est}, for every $\delta <2\alpha+1$,
\begin{equation}\label{eq-popotter}
H_2^\ep(|t-s|)=\frac {  H_2  (\ep |t-s|)}{  H_2  (\ep)}\le A_\delta\,|t-s|^{2\alpha+1-\delta}
\end{equation}
which gives
\begin{equation}\label{eq-maj1}
 \E(|Y^\ep_t-Y^\ep_s|^2)
%\le 2\, \frac{  H_2  (\ep|t-s|)}{  H_2  (\ep)}
\le 2 A_\delta\,|t-s|^{2\alpha+1-\delta}
\end{equation}
and, as $\delta$ can be chosen so that $2\alpha+1-\delta>0$, by
Kolmogorov's  compactness criterion the family $(Y^\ep)_\ep$ is
tight as $\ep\to 0$. 

Let us now consider the case $\alpha > 0 $.
Again, as for small $\ep$ we are led to the consideration of $h$ in
a right neighborhood of $0$ where it is increasing. This implies,
thanks to Theorem \ref{th:mod-con} iii),
$$
\E(|Y^\ep_t-Y^\ep_s|^2)\le H_2^\ep (t)-  H_2^\ep( s) +  H_2^\ep(|t-s|).
$$
Again Potter's bounds give \eqref{eq-popotter} for every $\delta <2\alpha+1$
%$$
% H_2^\ep(|t-s|)=\frac {  H_2  (\ep u)}{  H_2  (\ep)}\le A_\delta\,|t-s|^{2\alpha+1-\delta}
%$$
and by Lagrange's Theorem, for $s\leq t$
$$
H_2^\ep (t)-  H_2^\ep( s)= (  H_2^\ep)'(\xi_{s,t}) (t-s)=
\frac{\ep ( H_2 )' (\ep \xi_{s,t})}{  H_2  (\ep)} (t-s),
$$
for some suitable  $\xi_{s,t}\in[s,t]$. Furthermore thanks to \eqref{eq-bounded} applied to the $RV_{2\alpha+1}$ function $H_2$,
$u\mapsto\frac{\ep  ( H_2 )' (\ep u)}{  H_2  (\ep)}$ is bounded uniformly in $\ep$ and therefore there exists $C>0$ such that
\begin{equation}\label{eq-maj2}
\E(|Y^\ep_t-Y^\ep_s|^2)\le C|t-s|+ \frac{  H_2  (\ep|t-s|)}{  H_2  (\ep)}\le C |t-s|+
 A_\delta\,|t-s|^{2\alpha+1-\delta}
\end{equation}
and again by the Kolmogorov compactness criterion the family
$(Y^\ep)_\ep$ is tight as $\ep\to 0$. Continuity of
$X$ follows from \eqref{eq-maj1} and \eqref{eq-maj2}.

In order to conclude we just need to compute the limit of the covariance functions of the
processes $(Y^\ep)_\ep$. This follows from the next Proposition, that holds also for $\alpha=0$.
\cvd
\end{proof}
\begin{proposition}\label{prop-cinquequattordici} Let $Y^\ep$
be as in the
statement of Theorem \ref{teo-eta-1} for $\alpha>-\frac 12$ and let us denote
$k^\ep$ the covariance function of $Y^\ep$. Then 
\begin{equation}\label{eq-kep}
\lim_{\ep\to 0} k^\ep(s,t)=(2\alpha+1)\int_0^{t\wedge s}(t-u)^\alpha
(s-u)^\alpha\, du
\end{equation}
\end{proposition}
\begin{proof} We have for $s\leq t$ (recall that $k$ is the covariance function of $X$)
%$$
%k^\ep(s,t)=\frac 1{  H_2  (\ep)}\,k(\ep t,\ep s)
%$$
\begin{equation}\label{eq-kaep}
\begin{array}{c}
\displaystyle k^\ep(s,t)=\frac 1{  H_2  (\ep)}\,k(\ep t,\ep s)=\\[10pt]
\displaystyle =\int_0^{\ep s}\frac 1{  H_2  (\ep)}\,\sqrt{ h^2 (\ep t-u)}\sqrt{ h^2 (\ep s-u)}\,du=\\[12pt]
\displaystyle =\int_0^{ s}\sqrt{\frac{\ep\, h^2 (\ep (t-u))}{  H_2  (\ep)}}\sqrt{\frac{\ep\, h^2 (\ep (s-u))}{  H_2  (\ep)}}\,du.
\end{array}
\end{equation}
As $H_2\in RV_{2\alpha+1}$, \eqref{eq-deriv0} applied to $\lambda=  H_2  $  gives
$$
\frac{\ep\, h^2(\ep  t)}{ H_2  (\ep)}=\frac {\ep H_2'(\ep t)}{H_2(\ep)} \enspace\mathop{\to}_{\ep\to 0}(2\alpha+1)\, t^{2\alpha}
$$
i.e.
$$
\sqrt{\frac{\ep h^2 (\ep (t-u))}{  H_2  (\ep)}}\enspace\mathop{\to}_{\ep\to 0}\enspace \sqrt{2\alpha+1}\, (t-u)^{\alpha}.
$$
In order to take the limit in the integral some estimates are required: thanks to \eqref{eq-deriv} applied to $\lambda=H_2$,
in a right neighborhood of the origin we have
$$
\displaylines{
\frac{\ep\, h^2 (\ep(t-u))}{H_2  (\ep)}=\frac{\ep H_2'(\ep(t-u))}{H_2(\ep)}%=\cr
=\frac{H_2(\ep(t-u))}{H_2(\ep)}
\Bigl(\ep\frac {c'(\ep (t-u))}{c(\ep (t-u))}+\bigr(2\alpha+1-\rho(\ep(t-u))\bigr)\,\frac 1{t-u}\Bigr)
%\le\cr \le K_\delta\,\bigr((t-u)^{2\alpha+1 -\delta}+o((t-u)^{2\alpha -\delta})).\cr
}
$$
By Potter's bounds for every $\delta>0$ we have, for $0<\ep<T_\delta$ and some constant $C_\delta$
$$
\frac {  H_2  (\ep (t-u))}{  H_2  (\ep)}\le C_\delta (t-u)^{2\alpha+1-\delta}
$$
so that  for small $\ep$
$$
\displaylines{
\frac{\ep\, h^2 (\ep(s-u))}{H_2  (\ep)}\le
C_\delta\Bigl(\ep\frac {c'(\ep (s-u))}{c(\ep (s-u))}+\bigr(2\alpha+1-\rho(\ep(s-u))\bigr)\,\frac 1{s-u}\Bigr)(s-u)^{2\alpha+1-\delta}\le\cr \le K_\delta\,\bigr((s-u)^{2\alpha -\delta}+o((s-u)^{2\alpha -\delta}))\le (K_\delta+1)(s-u)^{2\alpha -\delta}.\cr
}
$$

As  we can choose $\delta$ so that $\alpha-\frac12\,\delta>-1$, we can apply Lebesgue's
theorem and pass to the limit in \eqref{eq-kaep} and obtain \eqref{eq-kep}.
\cvd
\end{proof}
The following statement takes care of the case $\alpha=0$.
\begin{theorem}\label{teo-eta-2}
Let us assume $h\in RV_0$ and that Assumption \ref{as-1} is satisfied. Assume moreover that at least one of the following conditions is satisfied in a right neighborhood of the origin.
\tin{a)} $\rho>0$;
\tin{b)} $\rho<0$;
\tin {c)} $\rho=0$ and $c'\le 0$;
\tin {d)} $\rho=0$ and $c'\ge 0$.

If $Y^\ep_t=\frac 1{\sqrt{H_2(\ep)}}\,X_{t\ep}$ then the family $(Y^\ep)_\ep$ converges to  a Brownian motion.
\end{theorem}
\begin{proof} As in the proof of Theorem \ref{teo-eta-2} we have
$$
 Y^\ep_t=\int_0^{t} h_\ep(t-u)\, dB_u,
$$
where
$$
h_\ep(t)=\frac 1{\sqrt{H_2(\ep)}}\, \sqrt{\ep}\,h(\ep t)
$$
so that
$$
H^\ep_2(t):=\int_0^t h_\ep^2(u)\, du=\frac 1{H_2(\ep)}\int_0^{\ep t}
\ep h^2(\ep u)\, du=\frac 1{H_2(\ep)}\int_0^{\ep t} h^2(u)\,
du=\frac {H_2 (\ep t)}{H_2 (\ep)}\cdotp
$$
First let us consider cases a) and c). As for small values of $\ep$ the behavior of $h^2_\ep$ reduces to
the consideration of $h^2$ in a right neighborhood
of $0$, by Lemma \ref{lem-zero} a) $h_\ep^2$ is non increasing and this implies,
by Theorem \ref{th:mod-con} ii), that
$$
\E(|Y^\ep_t-Y^\ep_s|^2)\le 2H_\ep^2(|t-s|).
$$
By Potter's bounds, \eqref{Potter-est}, for every $\delta <1$,
$$
 H_2^\ep(|t-s|)=\frac {  H_2  (\ep |t-s|)}{  H_2  (\ep)}\le A_\delta\,|t-s|^{1-\delta}
$$
which gives
$$
\E(|Y^\ep_t-Y^\ep_s|^2)\le 2\, \frac{  H_2  (\ep|t-s|)}{  H_2  (\ep)}\le 2 A_\delta\,|t-s|^{1-\delta}
$$
and, as $\delta$ can be chosen so that $1-\delta>0$, by Kolmogorov's criterion the family $(Y^\ep)_\ep$ is tight as $\ep\to 0$.

Let us now consider cases b) and d), so that $h^2$ is non decreasing in a right neighborhood of $0$.
Again as for small $\ep$ we are led to the consideration of $h$ in a right neighborhood of $0$, this implies that, thanks to Theorem \ref{th:mod-con},
$$
\E(|Y^\ep_t-Y^\ep_s|^2)\le H_2^\ep (t)-  H_2^\ep( s) +  H_2^\ep(|t-s|).
$$
By Potter's bounds, \eqref{Potter-est}, for every $\delta <1$,
$$
 H_2^\ep(|t-s|)=\frac {  H_2  (\ep |t-s|)}{  H_2  (\ep)}\le A_\delta\,|t-s|^{1-\delta}
$$
whereas by Lagrange's Theorem, for $s\leq t$,
$$
 H_2^\ep (t)-  H_2^\ep( s)= (  H_2^\ep)'(\xi_{s,t}) (t-s)=
\frac{\ep ( H_2 )' (\ep \xi_{s,t})}{  H_2  (\ep)} (t-s),
$$
for some suitable  $\xi_{s,t}\in[s,t]$. Furthermore thanks to \eqref{eq-bounded} applied to the $RV_{1}$ function $H_2$,
$\frac{\ep  ( H_2 )' (\ep \xi_{s,t})}{  H_2  (\ep)}$ is bounded and therefore there exists $C>0$ such that
$$
\E(|Y^\ep_t-Y^\ep_s|^2)\le 2\, \frac{  H_2  (\ep|t-s|)}{  H_2  (\ep)}\le A_\delta\,|t-s|^{1-\delta} + C |t-s|
$$
and again by Kolmogorov's criterion the family $(Y^\ep)_\ep$ is tight as $\ep\to 0$.

Proposition \ref{prop-cinquequattordici} allows to compute the limit of the
covariance functions of $Y^\ep$ and to conclude.
\cvd
\end{proof}
\begin{cor}\label{cor-st-rough}
In the hypotheses of Theorem \ref{teo-eta-1} or Theorem \ref{teo-eta-2},
the family $(X^\ep)_\ep$, where $X_t^\ep=X_{t\ep}$ for $t\in[0,1]$, enjoys an LDP at speed $g(\ep)=\frac 1{H_2(\ep)}$ and with rate function defined in \eqref{eq-rkhs0},
$\cl Y^0$ denoting the RKHS  of the process
$$
Y^0_t=\sqrt{2\alpha+1} \int_0^t(t-u)^{\alpha}\, dB_u.
$$
(which is a Brownian motion if $\alpha=0$).
\end{cor}
\begin{proof} Just note that we can write
$$
X_{\ep t}=\sqrt{H_2(\ep)}\Bigl(\frac 1{\sqrt{H_2(\ep)}}\, X_{\ep t}\Bigr)=\sqrt{H_2(\ep)}Y^\ep_t
$$
and apply Theorems \ref{teo-eta-1} or \ref{teo-eta-2}  and Corollary \ref{cor-st}.
\cvd
\end{proof}
Note that the rather ``botanical'' statement of Theorem \ref{teo-eta-2} is intended to
ensure that the integrand $h$ is monotonic near $0$ and, more or less, our results hold under this assumption that is required in order to obtain tightness through Theorem \ref{th:mod-con}.

In particular we are not able to treat integrators $h$ that are highly oscillating near $0$, as $t\mapsto 1+t\sin \frac1t$, which is a $RV_0$ function with $\rho=0$, but with an unbounded $c'$.

Note also that the limit process $ Y^0$ in Corollary \ref{cor-st-rough} is self similar, in agreement with Corollary \ref{cor-ss}. 
\begin{example}\label{exam-fin} \rm a) Let
\begin{equation}\label{eq-exam-rv}
h(t)=t^\alpha\,(-\log t)^\beta
\end{equation}
for some $\alpha>-\frac 12$. In the representation \eqref{eq-repr}
we can choose $t_0=1/e$, and $\rho=\frac{\beta}{-\log t}$. Therefore
the family of processes $(X^\ep)_\ep$  enjoys an LDP with rate
function given by the square of the RKHS of the
process $Y^0$ of Corollary \ref{cor-st-rough} divided by $2$ (which, note, depends on $\alpha$ only and is a self similar process).

The speed function is
$$
g(\ep)\sim \frac 1{H_2(\ep)}
$$
with
$$
H_2(\ep)=\int_0^\ep u^{2\alpha} (-\log u)^{2\beta}\, du.
$$
Note that, integrating by parts,
$$
H_2(\ep)=\frac 1{2\alpha+1}\,\ep^{2\alpha+1}(-\log \ep)^{2\beta}+\frac {2\beta}{(2\alpha+1)}\int_0^\ep u^{2\alpha}(-\log u)^{\beta-1}\, du\sim \frac 1{2\alpha+1}\,\ep^{2\alpha+1}(-\log \ep)^{2\beta}.
$$
We find here the same phenomenon as noted in Example
\ref{rem-conjectured-speed}: for different values of $\beta$ the
process $X$ exhibits an LDP behavior in small time with a common
rate function but different speeds.
\tin{b)} $h(t)=(\frac 1t\,\sin t)^{\gamma}$. Now $h\in RV_0$ with
the simple representation $c(t)=h(t)$ and $\rho=0$. As $\frac
1t\,\sin t=1-\frac 16\, t^2+o(t^2)$, near $0$ $h$ is increasing for
$\gamma<0$ and decreasing for $\gamma>0$, so that conditions b) or
d) of Theorem \ref{teo-eta-2} are satisfied and the family of small
time processes $(X^{\ep })_\ep$ enjoys an LDP with the rate function
as in \eqref{eq-rkhs0} for a Brownian motion.  Also the speed here
is independent of $\gamma$, as
$$
g(\ep)=\Bigl(\int_0^\ep\bigl(\frac 1u\, \sin u\bigr)^\gamma\,du \Bigr)^{-1}\sim \frac 1\ep\cdotp
$$
\end{example}
In Example \ref{exam-fin} a) one might consider the integrand $h$ in \eqref{eq-exam-rv}
with $\alpha=-\frac 12$ which, for $\beta>1$ gives rise to a square integrable
function so that the process $X$ in \eqref{eq-rv0} is well defined. In the next
section we see that in this case the LDP behavior of $X$ in small time is
much different.
\subsection{The super rough case: $h\in RV_\alpha$, $\alpha=-\frac 12$}
\label{par-superrough}
In this section  we investigate the small time LD asymptotics of a
Volterra process $X$ as in \eqref{eq-rv0} when $h\in RV_{-1/2}$ and
make a connection with the super rough volatility models as
introduced in Gulisashvili \cite{Gu-SR}, that we now recall. Note that in this case $H_2\in RV_0$.

A modulus of continuity is an {\it increasing} function
$\eta:\R^+\to\R^+$, vanishing at $0$ and continuous at $0$. Let $X$
be a real square integrable process. A modulus of continuity $\eta$
is a {\it weak modulus of continuity} for $X$ if for every $s,t$
$$
\E(|X_s-X_t|^2)\le \eta^2(|t-s|)
$$
It is a result of Fernique \cite{Fer} (see also Dudley \cite{Dud})
that if
\begin{equation}\label{eq-cond-fernique}
\int_*^{+\infty}\eta(\e^{-x^2})\, dx
=\int_0^* \eta(u)(-\log u)^{-1/2}\, \frac 1u\, du<+\infty
\end{equation}
then $X$ has a continuous version.

If $X$ is a Volterra process as in \eqref{eq-rv0} then, by Theorem
\ref{th:mod-con} ii) if $h^2$ is non increasing,  $H_2$ is a weak
modulus of continuity for $X$. Note that then $H_2$ is a concave
function and therefore subadditive; so we have, thanks to Theorem
\ref{th:mod-con} i), for $t\ge s$,
$$
\E(|X_s-X_t|^2)\ge H_2(|t-s|)\ge H_2(t)-H_2(s)
$$
and by Lemma 7.2.13 of \cite{marcus-rosen} we obtain  the lower
bound for the modulus of continuity of the paths of $X$
$$
\lim_{\delta\to 0}\sup_{u\le\delta}\frac{|X_u|}{\sqrt{H_2(u)}\log\log (H_2(u)^{-1/2})}\ge\sqrt{2}.
$$
In particular if $H_2$ grows at $0$ faster than $t\mapsto t^\alpha$
for any power $\alpha>0$, then $X$ cannot have H\"older continuous
paths of any exponent. Processes with paths enjoying this kind of
lack of regularity were introduced in Gulisashvili \cite{Gu-SR} and
are called {\it super rough}. They are intended to
provide models for the volatility in the evolution of prices in
financial models as sometimes data suggest severe lack of
regularity.
\begin{example}\label{exam-u1}\rm a) Let $X$ be as in \eqref{eq-rv0} with
$h(t)=t^{-1/2}(-\log t)^{-\gamma}$. If $\gamma>\frac 12$, then $h$ is square integrable near $0$ and $h^2(t)=\frac 1t\,(-\log t)^{-2\gamma}$ is non increasing near $0$, hence, thanks to Theorem \ref{th:mod-con} ii)
$$
\E(|X_s-X_t|^2)\le 2H_2(|t-s|)
$$
where
$$
H_2(t)=\int_0^t h^2(u)\, du=\int_0^t\frac 1u\,(-\log u)^{-2\gamma}\,
du= \frac 1{2\gamma-1}\,(-\log t)^{-2\gamma+1}.
$$
$X$ has therefore weak modulus of continuity
$\eta=\sqrt{2H_2}=(2\gamma-1)^{-1/2}\,(-\log t)^{\frac 12-\gamma}$.
By the above mentioned Fernique's result $X$ has a continuous
version as soon as $\gamma>1$. Note that $\eta\in RV_0$.
\tin{b)} Conversely a super rough Gaussian process can  be obtained
starting from a modulus of continuity. Actually, let $\eta$ be a
modulus of continuity satisfying Fernique's condition
\eqref{eq-cond-fernique}, then, if $h(t)=\sqrt{(\eta^2)'(t)}$ is
square integrable, the process $X$ as in \eqref{eq-rv0} has $\eta$
as its weak modulus of continuity. Moreover the additional
assumption of concavity of $\eta^2$ will ensure that the modulus of
continuity of the paths of $X$ grows at $0$ faster than $\eta$. Note
that in this situation, if $\eta\in RV_0$, then $h\in RV_{-1/2}$.
Let for instance
$$
\eta(t)=(-\log x)^{-\beta} \, (\log(-\log x))^{\alpha},
$$
for $\beta>0$ and $\alpha\in \R$.
It is straightforward to show that $\eta^2$ is a concave function and that $\eta$ satisfy
the Fernique condition if $\beta>\frac12$.
\end{example}
Let us go back to the investigation of the LD behavior of these
super rough models in small time. Let $h$ be a square integrable
integrand such that $t\mapsto h^2(t)$ is non increasing, so that
$\eta^2:=H_2$ is a concave function. Notice that in this case
$\eta^2(t)=k(t,t)=H_2(t)$. Recall that $H_2$ must belong to  $ RV_0$, i.e.
\begin{equation}\label{eq:slo-var}
\lim_{\ep\to 0} \frac{H_2(\ep t)}{H_2(\ep)}=1\ .
\end{equation}

We prove now that if a process $X$ has a weak modulus of continuity $\eta$ belonging to $RV_0$, then it cannot enjoy an LDP in small time. It enjoys, however, a weak LDP (i.e. with upper bound for compact sets only).

Let $X$ be as in \eqref{eq-rv0}.   The first step is to investigate the convergence of the family of processes $(Y^\ep)_\ep$ where, as above,
$$
Y^\ep_t=\frac 1{\sqrt{H_2(\ep)}}\, X_{\ep t}=\frac 1{\eta(\ep)}\, X_{\ep t}.
$$ 
We have, for $s<t$, ($k=$the covariance function of  $X$)
$$
\begin{array}{c}
\displaystyle k(t,s)=\int_0^{s} h(  t-u) h( s-u)\, du
\leq  \Bigl(\int_0^{s} h^2( t-u)\, du \Bigr)^{\frac 12 }\Big(\int_0^{s} h^2(  s-u)\,du\Big)^{\frac 12 }=\\ [10pt]
=\bigl(H_2(t)-H_2(t-s)\bigr)^{\frac 12 }\sqrt{H_2(s)}.
\end{array}
$$

\begin{proposition}\label{lemma:cov-lim} Let $X$ be as in \eqref{eq-rv0} with
an integrand $h\in RV_{-1/2}$ which is square
integrable and such that $\eta=\sqrt{H_2}$ satisfies Fernique's
condition \eqref{eq-cond-fernique}. We have,
$$
k_0(t,s):=\lim_{\ep\to 0}\frac 1{H_2(\ep)}\, k(\ep t, \ep s)=
\begin{cases}
1&\mbox{if }s=t\\
0&\mbox{if }s\not=t\ .
\end{cases}
$$
\end{proposition}

\begin{proof} Let $ k^{\ep}$ be the  covariance function of $Y_\ep$, i.e.
$$
k^{\ep}(t,s):=\frac 1{H_2(\ep)}\, k(\ep t, \ep s).
$$
For $s=t$ the result follows from \eqref{eq:slo-var},
as  $ k^\ep(t,t)=H_2(\ep t)$ and, as $H_2$ is
slowly varying. For $s\not=t$ we have
$$
k^{\ep}(s,t)\le \frac{\sqrt{H_2(\ep s)}}{\sqrt{H_2(\ep)}}\Bigl(\frac{H_2(\ep t)}{H_2(\ep)}-\frac{H_2(\ep (t-s))}{H_2(\ep)}\Bigr)^{\frac 12 }
$$
and, again as $H_2\in RV_0$, we have
$$
\lim_{\ep\to 0}  k^{\ep}(s,t)=0.
$$
\scvd
\end{proof}
As the limit $k_0$ is not the covariance of a continuous Gaussian
processes, the family of processes $(Y^\ep)_\ep$ does not converge
as $\ep\to 0$ in the space of continuous functions and the results
of  \S\ref{par-smalltime} do not apply in order to obtain a pathwise
LDP for the small time family $(X^{\ep})_\ep$.

Actually the following more precise negative result holds.
\begin{theorem}\label{th-negative}Let $X$ be as in \eqref{eq-rv0}
with an integrand $h\in RV_{-1/2}$ which is square integrable and
such that $\eta=\sqrt{H_2}$ satisfies Fernique's condition
\eqref{eq-cond-fernique}. Then the family of small time processes
$(X^\ep)_\ep $ is not exponentially tight and does not enjoy an LDP
as $\ep\to 0$ at speed $g=H_2$.
\end{theorem}
\begin{proof}
Let for $n\ge 1$, $0\leq t_1< t_2<\dots t_n\leq 1$ and let us
consider the family of $n$-dimensional r.v.'s
$$
X_n^{\ep}=\bigl(X(\ep t_1),\ldots,X(\ep t_n)\bigr).
$$
As by Proposition \ref{lemma:cov-lim} the covariance matrix of
$(H_2(\ep)^{-1/2}X_n^\ep)_\ep $ converges to the identity matrix
and, by classical results on LD's for Gaussian vectors (the
Ellis-G\"artner theorem e.g.) this family satisfies an LDP, on
$\R^n$, with speed $H_2(\ep)^{-1}$ and rate function, for any
$t_1,\ldots,t_n$, given by,
\begin{equation}\label{eq:RF-fin-dim}
I_{n}(x_1,\ldots,x_n)= \frac 12 \sum_{i=1}^n x_i^2 .
\end{equation}
{\it If} $(X^\ep)_\ep $ enjoyed an LDP at speed $g=H_2^{-1}$ as
$\ep\to 0$, then by the contraction principle $(X_n^\ep)_\ep $ would
enjoy an LDP at speed $g$ with respect to a rate function $I$ such
that
\begin{equation}\label{eq-pseudo-I}
I_n(x_1,\dots,x_n)=\inf\{I(\gamma);
\gamma(t_1)=x_1,\dots,\gamma(t_n)=x_n\}.
\end{equation}
This together with \eqref{eq:RF-fin-dim} gives
$$
I(\gamma)\ge \frac 12 \sum_{i=1}^n \gamma(t_i)^2
$$
and, as this holds for every $n\ge1$, $x_1,\dots,x_n$, we have
\begin{equation}\label{eq-ideg}
I(\gamma)=\begin{cases}
0&\mbox{if }\gamma=0\\
+\infty&\mbox{otherwise}.
\end{cases}
\end{equation}
But this is absurd, as \eqref{eq-pseudo-I} would give now
$$
I_n(x_1,\dots,x_n)=\begin{cases}0&\mbox{if } x_1=\dots=x_n=0\cr
+\infty&\mbox{otherwise}
\end{cases}
$$
in contradiction with \eqref{eq:RF-fin-dim}. Therefore the family
$(X^\ep)_\ep$ does not enjoy an LDP at speed $H_2^{-1}$ and cannot
be exponentially tight.
\end{proof}
\begin{remark}\rm Theorem \ref{th-negative} states that the process $X$ cannot enjoy an LDP in small time at the natural speed $g=H_2^{-1}$. What about an LDP at a different speed?

At a speed $g$ that is faster than $H_2^{-1}$ it is easy to see that $(X^\ep_n)_n$ cannot enjoy an LDP: the Ellis G\"artner theorem would give $I_n\equiv 0$ and by \eqref{eq-pseudo-I} also $I\equiv 0$, which is not a rate function.

At a speed $g$ that is slower than $H_2^{-1}$ it turns out that 
$(X^\ep_n)_n$ enjoy an LDP with respect to a degenerate rate function, so that, if $(X^\ep)_\ep$ enjoyed an LDP in small time, then it with be with respect to a degenerate rate function, which is of limited interest in our opinion.

See the proof of Proposition \ref{th:diff-speed} for details
\end{remark}

\begin{remark}\rm Let us sum up the results of this section concerning
the super rough setting.

We have proved that no LDP holds for the process defined in
\eqref{eq-rv0} and with  respect to an integrand  $h$ such that
$H_2$ is slowly varying.

However, by the Dawson-G\"artner Theorem (see Theorem 4.6.1. in
\cite{DZ:97}), we have a pathwise LDP at speed $H_2^{-1}$ and rate
function $\widehat I$, defined in equation (4.6.2) in \cite{DZ:97},
in the space of real function on $[0,1]$ endowed with the topology
of pointwise convergence. Note that  for every continuous function
$\widehat I(\gamma)=I(\gamma)$ where $I$ is as in \eqref{eq-ideg}.

As proved above no LDP holds for the process defined in \eqref{eq-rv0} and with respect to an integrand $h$ such that $H_2$ is slowly varying.
% as we are not able to prove the exponential tightness of the family $(X^\ep)_\ep$ with $X^\ep_t=X_{\ep t}$.
We can however obtain a weaker result, i.e. a lower bound for open
sets and an upper bound for compact sets in the space of continuous
functions (i.e. a weak LDP).

Actually, let $G\subset \cl C$ be an open set in the topology of uniform convergence. If $0\in G$ then
$$
\liminf_{\ep\to 0} H_2(\ep) \log \P (X^\ep\in G)=0.
$$
If  $0\notin G$ instead, then the inequality
$$
\liminf_{\ep\to 0} H_2(\ep) \log \P (X^\ep\in G)\geq -\infty
$$
is obvious so that in any case we have the lower bound
$$
\liminf_{\ep\to 0} H_2(\ep) \log \P (X^\ep\in G)\geq -\inf_{\gamma\in G} I(\gamma).
$$
If $K\subset \cl C$ is a compact set in the topology of uniform convergence, then $K$ is closed  in the topology of pointwise convergence and  therefore the upper bound for compact sets
$$
\limsup_{\ep\to 0} H_2(\ep) \log \P (X^\ep\in K)\leq -\inf_{\gamma\in K} \hat I(\gamma)=-\inf_{\gamma\in K} I(\gamma)
$$
follows from the upper bound in the topology of the pointwise convergence mentioned above.
\end{remark}

\paragraph{Funding.}
The authors are supported by MIUR Excellence Department Project awarded to the Department of Mathematics, University
of Rome Tor Vergata (CUP E83C18000100006), by University of Rome Tor Vergata (project "Asymptotic Methods in 
Probability" (CUP E89C2000068\-00\-05) and project "Asymptotic Properties in Probability" (CUP E83C22001780005)) and 
by Indam-GNAMPA.
%%%%%%%%%%%%%%%%%%%%%%%%%%%%%%%%%%%%%%%%%%%%%%%%%%%%%%%%%

\bibliography{bibbase}
\bibliographystyle{amsplain}

\end{document}